\documentclass[12pt]{amsart}
\usepackage{amsmath,amscd,latexsym,verbatim,amssymb}
\usepackage{times}       

 \newcommand{\resp}{{\it resp.} }
\newcommand{\cf}{{\it cf.} }
\newcommand{\ie}{{\it i.e.} }

\newcommand{\Q}{\mathbf{Q}}

\newcommand{\C}{\mathbf{C}}
  \newcommand{\Z}{\mathbf{Z}}

\newcommand{\sD}{{\mathcal{D}}}
\newcommand{\sE}{{\mathcal{E}}}

\newcommand{\sH}{{\mathcal{H}}}
\newcommand{\sI}{{\mathcal{I}}}

\newcommand{\sO}{{\mathcal{O}}}

\newcommand{\sY}{{\mathcal{Y}}}

\newcommand{\inj}{\hookrightarrow}

  \newcommand{\ord}{{\rm{ord}}}

\renewcommand{\epsilon}{\varepsilon}
\renewcommand{\phi}{\varphi}

\renewcommand{\lim}{\varprojlim}

\font\sm=cmr10 at 10pt
\font\smit= cmti10 at 10pt
 
\newcommand{\Ker}{\operatorname{Ker}}

\newcommand{\Spec}{\operatorname{Spec}}

\newcounter{spec}
{\end{list}}

\swapnumbers

\newtheorem{thm}{Theorem}[subsection]
\newtheorem{lemma}[thm]{Lemma}
\newtheorem{prop}[thm]{Proposition}
\newtheorem{cor}[thm]{Corollary}

\theoremstyle{definition}

\newtheorem{ex}[thm]{Example}

\newtheorem{rem}[thm]{Remark}

\numberwithin{equation}{section}

\font\sm=cmr10 at10pt
 at12pt

\newcommand{\prf}{\noindent {\bf Proof. }}
\renewcommand{\qed}{\hfill $\square$\medskip}

\begin{document}

\title{An algebraic proof of Deligne's regularity criterion for integrable connections}
\author{Yves
Andr\'e}
  \address{D\'epartement de Math\'ematiques et Applications, \'Ecole Normale Sup\'erieure  \\ 
45 rue d'Ulm,  75230
  Paris Cedex 05\\France.}
\email{andre@dma.ens.fr}
  \date{today}
 \begin{abstract} Deligne's regularity criterion for an integrable 
connection $\nabla$ on a smooth complex algebraic variety $X$ says 
that $\nabla$ is regular along the irreducible divisors at infinity 
in some fixed normal compactification of $X$ if and only if the 
restriction of $\nabla$ to every smooth curve on $X$ is fuchsian ({\it 
i. e.} has only regular singularities at infinity). The ``only if" 
part is the difficult implication. Deligne's proof is transcendental  
and uses Hironaka's resolution of singularities. 

Following \cite{A}, we present a purely algebraic proof of this implication which does not use resolution  beyond the case of plane curves. It relies upon a study of the formal structure of integrable connections on surfaces with (possibly irregular) singularities along a divisor with normal crossings. 
  \end{abstract}
\maketitle



\bigskip

\begin{sloppypar}

\section{Introduction: fuchsian connections.}\label{s1}
 
\medskip

\subsection{} Let $X$ be a connected algebraic complex manifold, and let ${\sE} $ be an algebraic vector bundle on $X$ endowed with an integrable connection $\nabla$.

When $X$ is a curve, the dichotomy between regular and irregular 
singularities (at infinity) goes back to the 19th century. The 
connection $\nabla$ is said to be {\it fuchsian} if all singularities 
(at infinity) are regular, see Manin's classical paper \cite{M}. This is 
obviously a birational notion.

\medskip In higher dimension, one may consider a normal 
compactification $\bar X$ of $X$ and look at  the 
irreducible components $Z_j$ of $\partial \bar X = \bar X \setminus X$ of 
codimension one in $\bar X$. The local subring ${{\sO}_{X,Z_j}}$ of $\C(X)$ is a discrete valuation ring with residue field ${{\sO}_{X,Z_j}}/{{\frak m}_{X,Z_j}}=\C(Z_j)$. One is then in the familiar one-dimensional situation, over $\C(Z_j)$ instead of $\C$, and the notion of regularity of $\nabla$ along
$Z_j$ is defined in the usual way. Namely, consider a germ of vector field  $\theta_j$  on $X$ - viewed as a derivation of $\C(X)$ - which preserves ${\frak m}_{X,Z_j}$, but does not send it into ${\frak m}^2_{X,Z_j}$,  and which acts trivially on the residue field. Then $\nabla$ is said to be {\it regular along $Z_j$} if there is an ${{\sO}_{X,Z_j}}$-lattice in the generic fiber $\sE_{\C(X)}$ which is stable under $\nabla(\theta_j)$ (this condition does not depend on $\theta_j$,  cf. \cite[I.3.3.4]{AB}).

 \smallskip In order to obtain a birational notion of fuchsianity, one is then led 
to say that $\nabla$ is {\it fuchsian} if for any $(\bar X, Z_j)$ as before, 
$\nabla$ is regular at $Z_j$.
  
  At first look, this definition is rather forbidding, since it 
requires to consider all divisorial valuations of $\C(X)$ at the same 
time.
   Fortunately, it turns out that it suffices to consider only one normal 
compactification $\bar X$. Indeed, according to P. Deligne, one has the following 
characterizations of fuchsianity.

\begin{thm}\label{th1}  \cite[II.4.4, 4.6]{D} The following are equivalent:

\smallskip $i)$ $\nabla$ is fuchsian,

\smallskip $ii)$ for some normal compactification $\bar X$, 
$\nabla$ is regular along all irreducible 
components of $\partial \bar X$ of codimension one in $\bar X$,

\smallskip $iii)$ for any smooth curve $C$ and any locally 
closed embedding  $h: C\to X$, $h^\ast\nabla$ is fuchsian,

\smallskip $iv)$ for any smooth $Y$ and any morphism $f: Y\to 
X$, $f^\ast\nabla$ is fuchsian,

\smallskip $v)$ for some dominant morphism $f: Y\to X$ with 
$Y$ smooth, $f^\ast\nabla$ is fuchsian.
\end{thm}

Note that we {\it do not} assume that $\partial \bar X$ has normal crossings.
The difficult implication is $ii)\Rightarrow iii)$. The implication 
$iii)\Rightarrow i)$ is comparatively easy to establish ({\it cf. 
e.g.} \cite[I. 3.4.7]{AB}), and the other implications follow very easily 
from these two.
  The difficulty with $ii)\Rightarrow iii)$ arises when the closure of 
$C$ in $\bar X$ does not meet $\partial \bar X$ transversally. A closely 
related difficulty, with $ii)\Rightarrow i)$, is to show that 
$\nabla$ remains regular at the exceptional divisor when one blows up 
a subvariety of $\partial \bar X$.

\subsection{} Deligne's proof (as contained in the erratum to \cite{D}) is not elementary: it relies upon certain transcendental complex-analytic arguments on one hand, and upon Hironaka's resolution of singularities on the other hand. 

\smallskip Recently, inspired by part of Sabbah's work on asymptotic analysis in dimension 2 \cite{S}, the author has found a purely algebraic proof of Deligne's regularity criterion, as a consequence of a more general result about the semicontinuity of the Poincar\'e-Katz rank \cite{A}. 

The aim of this text is to explain this argument  in the simplified form needed for \ref{th1}, \ie in the context of algebraic connections which are regular at infinity. We do not use resolution of singularities beyond resolution of plane curves.

\medskip We refer to \cite{B} for a nice introduction to the background and the eventful story of this problem. Its difficulty is related to the fact that the standard techniques of logarithmic lattices break down in this algebraic context; along the way, it becomes necessary to go beyond the framework of regular connections.
 
\subsection{} More precisely, we present an algebraic proof of the following local refined form of Deligne's regularity criterion\footnote{this statement appears in \cite[I.5.4]{AB}, but the argument given there works only in case the polar divisor has normal crossings. Indeed, as was pointed out by J. Bernstein, lemma 5.5 of op. cit. on which it relies does not hold in greater generality.}.

\begin{thm}\label{t2}  Let $X $ be a normal connected
algebraic variety over $\C$. Let $U$ be 
a smooth open subset of $X$, with complement ${\partial}X = 
X  \setminus U$, and let $Q$ be a closed point of $\partial X$. Let $Z_1, \dots, Z_t$ be the irreducible components of $\partial X$ of codimension one in $X$ which pass through $Q$.
  \hfill\break Let $C$ be a smooth connected curve, $P\in C$ be a closed point, and $h: C \to  
X $  be a
morphism such that  $h( C \setminus P) \subset U$ and $ h(P)= Q $.
    \hfill\break Let  $ {\sE} $ be a vector bundle on $U$ with an integrable connection $\nabla$.
 If $({\sE}, \nabla)$ is regular along each $Z_j$, then the vector bundle 
$h^\ast({\sE},\nabla ) $ on $C\setminus P$ is regular at $P$.
\end{thm}

\subsection{Warning.}\label{w} What does regularity actually mean? {\it Quot capita tot sensus}: in the literature, one can find {\it analytic} definitions (moderate growth of solutions of meromorphic connections), {\it algebraic} definitions (such as the items in \ref{th1}, reducedness of the characteristic variety, and many others \cite{AB}), as well as {\it mixed} definitions (of GAGA type: formal/algebraic versus analytic properties of connections).

The occurrence of so many different definitions is an obvious sign of the richness of the concept, but at the same time an amazing source of confusion.  
All of them are supposed to be equivalent in their common domains of application, but this is often a matter of folklore or belief, as a complete and precise dictionary is still lacking (especially in the sensitive case where the polar divisor has non-normal crossings). A strange situation, indeed, for a supposedly well-understood classical notion! 

\smallskip It is thus essential to keep in mind the definition which we have adopted right at the beginning in order to understand what we do in this paper: comparing, in a purely algebraic way, a few of the existing algebraic definitions of {\it regular} connections. The wit is that the core of the paper deals with {\it irregular} connections.

 \section{Reduction to the plane case.}\label{s2}

\subsection{  } In this section, we reduce \ref{t2} to the case where $X $ is the projective plane\footnote{the next two lemmas are joint work wih F. Baldassarri.}.

We may assume $\dim X\geq 2$, otherwise the theorem is essentially trivial. Replacing $X$ by a quasi-projective neighborhood of $Q$ and taking the normalization of its Zariski closure in some projective embedding, we may assume that $X$ is projective.

\begin{lemma}\label{l1}  Let $X  \subset {\bf P}^N_\C$ be a closed normal connected subvariety of dimension $d \geq  2$.  Let $U$, $Z_j$, $Q$, $h:C\to X$, $P$ be as in \ref{t2}.
 
  For any sufficiently large integer $\delta$, there exists an irreducible
complete intersection $Y \subset {\bf P}^N_\C $ of dimension $N-d+2$ and multidegree $(\delta,
\dots,\delta)$ such that:
\hfill \break
$(i)$ $Y$
contains $h(C)$ (with
reduced induced structure) and cuts $U$ transversally at its generic point $\eta_{h(C)}$;
\hfill
\break
$(ii)$ $Y \cap X $ is an normal connected surface and $Y$
cuts $U\setminus (h(C)\setminus Q)$ transversally;
\hfill \break
$(iii)$ in a
neighborhood of $Q$, $Y$ cuts each $Z_j \setminus Q$
transversally, and does not cut any irreducible component of $\partial
X$ of codimension $> 1$ in $X $, nor the singular locus of
$\, Z := \cup  Z_j$, except in  $Q$.
 \end{lemma}

\prf This is more or less standard, but, for lack of adequate reference, we give
some detail. For short, we change a little notation and now write $C$ for the closure of $h(C)$ in $X$, with reduced structure (an integral, possibly singular, curve). 
  Let $\pi: {\tilde {\bf P}} \longrightarrow {\bf P}^N$ be the blow-up centered at $C $, and let $E\subset {\tilde {\bf P}}$ the exceptional divisor.  

    We set ${\sI}_{C } = {\Ker}
({\sO}_{{\bf P}^N} \longrightarrow {\sO}_{C }$).
    Then ${\rm Im}(\pi^\ast  {\sI}_{C}  \to {\sO}_{\tilde {\bf P}})= {\sO
 }(-E)$ and,   for  $\delta >>0$,
    $\pi^{\ast} ({\sO}_{{\bf P}^N}(\delta))
\otimes
{\sO}(-E)$ is very ample: a basis of sections
defines an embedding
of  $\tilde {\bf P}$ into  ${\bf P}^M$. On the other hand, since $\pi$ is birational and ${\bf P}^N$ is normal, one has $\pi_{\ast}{\sO}_{\tilde {\bf P}} =
{\sO}_{{\bf P}^N}$. It follows that the two natural arrows
$${\sI}_{C }  \to \pi_\ast({\rm Im}(\pi^\ast  {\sI}_{C}  \to {\sO}_{\tilde {\bf P}}))= \pi_\ast{\sO
 }(-E)\to  \pi_{\ast}{\sO}_{\tilde {\bf P}} =
{\sO}_{{\bf P}^N}$$ are inclusions of ideals of ${\sO}_{{\bf P}^N}$. Let $D\subset C$ be the closed subscheme corresponding to $\pi_\ast{\sO
 }(-E)$. If $D\neq C $,  $D$ would be punctual and $1\in {\sO}_{{\bf P}^N\setminus D^{red}}$ would correspond to a function on $ {\tilde {\bf P}}\setminus (\pi^{-1}(D))^{red}$ vanishing on $E\setminus (\pi^{-1}(D))^{red}$, a contradiction. Hence $D=C$ and therefore $\pi_{\ast} {\sO}(-E) = {\sI}_{C}$.  From the projection formula, one
deduces
$$\pi_{\ast}( \pi^{\ast}({\sO}_{{\bf P}^N}(\delta)) \otimes  {\sO}(-E))
\cong {\sO}_{{\bf P}^N}(\delta) \otimes {\sI}_{C}
\cong {\sI}_{C} (\delta)
,\;  $$
   whence
    $$   {\Ker}  (H^0({\bf P}^N, {\sO}_{{\bf P}^N}(\delta)) \rightarrow H^0(C, {\sO}_{C}(\delta))) =  
      H^0({\bf P}^N, \, {\sI}_{C}(\delta)) $$ $$ = H^0({\tilde {\bf P}}, \pi^{\ast}({\sO}_{{\bf
P}^N}(\delta)) \otimes {\sO}(-E))  ,   
    $$
and the linear system
of hypersurfaces of degree $\delta $ in ${\bf P}^N$ containing $C$
gives rise
     to a locally closed embedding
    $${\bf P}^N \setminus C
\hookrightarrow {\bf P}^M = {\bf P}( H^0({\bf P}^N, {\sI}_{C}(\delta)))
\;,$$
    with Zariski closure $\tilde {\bf P}$.
    The canonical bijection
between hyperplanes $\sH$ of ${\bf P}^M$ and
    hypersurfaces $H$
of degree $\delta$ in ${\bf P}^N$ containing $C$, is such that
the intersection
${\sH} \cap ({\bf P}^N \setminus C)$
(in  ${\bf P}^M$)
    equals $H \setminus C$.  So, the
intersection of $X \setminus C$ with a general
complete intersection $Y$ of multidegree $(\delta,\dots,\delta)$ ($1 \leq s
\leq d-1$ entries) in ${\bf P}^N$ containing $\eta_C$,  is the
intersection of $X \setminus C$ with a general
linear subvariety $\sY$ of codimension $s$ in ${\bf P}^M$. By
\cite[Exp. XI, Thm. 2.1. $(i)$]{sga},
$Y$ cuts $X \setminus C$
(resp.  the smooth part of  $Z \setminus (C \cap
Z)$) transversally  and  intersects properly  any
irreducible component of $\partial X \setminus (C \cap
\partial X)$.  Since $s < d$, Bertini's theorem shows that the
intersection of ${\sY}$ with the strict transform of
$X $ in ${\bf P}^M$ is normal and connected. On the other hand, since
$\eta_C$ is a simple point of $X$, it is well-known that a general
complete intersection of $s$ hypersurfaces of degree $\delta$ in ${\bf
P}^N$ containing $\eta_C$,  intersects $X$ transversally at this
point.
\par
  Applying these considerations for 
$s =d-2$, one obtains  $(i)$, $(ii)$
and $(iii)$.  \qed

\begin{cor}\label{c1}  There
exists a normal quasi-projective irreducible neighborhood $X'$ of
$Q$ in $Y \cap X $, containing an open subset of
$h(C) $, such that $U' :=
X' \cap U$ is smooth, the
distinct irreducible components  of codimension one of
${\partial}X' = X' \setminus U'$ passing through $Q$, are
precisely $ Z_1 \cap X', \dots,  Z_t \cap X'$. \qed
\end{cor}

If $({\sE}, \nabla)$ is regular along each $Z_1,\dots Z_t$, so is its pullback $({\sE}, \nabla)_{\mid U'}$ along $ Z_1 \cap X', \dots,  Z_t \cap X'$ (\cf \cite[I.3.4.4]{AB}). This reduces theorem  \ref{t2} to the case where $X $ is a normal surface, or even a projective normal surface 
(by the same argument as in the beginning of this section).

 \begin{lemma}\label{l2} In the notation of lemma \ref{l1}, let us further assume
that $d=2$. There exists a morphism
    $g: X \longrightarrow {\bf P}^2_\C$, which is finite in a neighborhood $V$ of $Q$, such that $g (h(C)\cap V) $ is not contained in the branch locus, and such that, for any irreducible component $T$
of $\partial X$ of dimension 1 with  $Q \notin T$,  $ g (Q) \notin
g(T)$. \end{lemma}

    \prf  Let ${\bf G}(N-3, {\bf P}^N)$ be the Grassmannian of 
linear subvarieties of ${\bf P}^N$ of codimension 3, and let $G$ be its dense open subset consisting of linear subvarieties which do not intersect $X$. 

For any $\gamma\in G$, $X$ may be considered as a closed subvariety of the blow-up $\tilde{\bf P}_\gamma$ of ${\bf P}^N$ at $\gamma$, and the projection with center $\gamma$  $$p_{//\gamma}:\, \tilde{\bf P}_\gamma \to {\bf P}^2$$ induces a morphism 
$$g_\gamma: X \to {\bf P}^2.$$ 

Let ${\bf \Lambda} \subset  {\bf G}(N-2, {\bf P}^N) \times {\bf G}(N-3, {\bf
P}^N) $  the incidence subvariety (locus of $(\lambda ,{\alpha})$
such that ${\lambda}$ contains $\alpha$), and let $p_2, p_3$ be the natural projections. Notice that  $p_3$ is a fibration with  fiber
${\bf P}^2$ and admits a section above $G$: there is a unique $\lambda_\gamma\in  {\bf G}(N-2, {\bf P}^N)$ of ${\bf P}^N$ of codimension 2  passing through $Q$ and containing $\gamma$. 
Then $\gamma$ varies, the $\lambda_\gamma$ form (via $p_2$) a dense open subset of ${\bf G}(N-2, {\bf P}^N)$.

On the other hand, $\lambda_\gamma$ may be identified with the fiber of $p_{//\gamma}$ above $g_\gamma(Q)\in {\bf P}^2 $, and $\lambda_\gamma \cap X$ with $g_\gamma^{-1}(g_\gamma(Q))$. 

By Bertini, one deduces that there is an open dense subset $G'\subset G$ such that for any $\gamma\in G'$, $g_\gamma$ is finite in a neighborhood of $Q$. 

 Moreover, if there exists $\lambda \in p_3^{-1}(\gamma)$ which intersects $X$ transversally and cuts $h(C)\cap V$ (\resp which passes through $Q$ and avoids the $T$'s), then $g_\gamma(h(C)\cap V) $ is not contained in the branch locus of $g_\gamma$ (\resp $g_\gamma(Q)\notin g_\gamma(T)$). 
 
 One deduces that there is an open dense subset $G''\subset G'$ such that for any $\gamma\in G''$, $g= g_\gamma$ satisfies the conditions in the lemma. \qed

 In the situation of theorem \ref{t2} with $X $ a projective
surface, $g$ induces an \'etale morphism 
$$ X' = V \setminus (g_{\mid V}^{-1}(B) \cup g_{\mid V}^{-1}(g (\partial X))))  
\longrightarrow
 {\bf P}^2_\C \setminus (B \cup g(\partial X)) $$ which is finite over its image. Up to replacing $C$ by a suitable neighborhood of $P$, $h(C\setminus P)\subset U\cap X'$. 

Moreover, using  \cite[I.3.2.5]{AB} and \ref{l2}, one sees that the push-forward $g_\ast(({\sE}, \nabla)_{\mid U\cap X'})$ on $g(U\cap X' )$ is regular along each 1-dimensional irreducible components of ${\bf P}^2_\C \setminus g(U\cap X' )$ passing through $g(Q)$ (which are either the $g(Z_j)$'s or else are contained in $B$).

On the other hand, if  the connection $(g\circ h)^\ast g_\ast(({\sE}, \nabla)_{\mid U\cap X'})$ on $C\setminus P$ is regular at $P$, the same is true for its subconnection $h^\ast(({\sE}, \nabla)_{\mid U\cap   X'})$. 

Therefore, in order to prove the theorem, one may substitute to $X$ any Zariski neighborhood of $ g(Q)$ in the projective plane. 

\begin{cor}\label{c2} Statement 1.3.1 holds in general if it holds when $X$ is an affine neighborhood of the origin $Q$ in the complex plane and  $\,\partial X=  Z := Z_1\cup \ldots \cup Z_t$.
\qed \end{cor}
 
\subsection{  }   By the classical embedded resolution of plane curves, there is a sequence of blow-ups
 $$\pi:  X'= X_N\to \cdots\to X_0= X$$
 such that $\pi^{-1}(h(C)\cup Z )$ (with its reduced structure) has strict normal crossings. We denote by $C', Z'_j$ the strict transforms of $h(C), Z_j$ respectively; they are smooth curves on $X'$, and $h$ lifts to a morphism $C\to C'$. 
  
\smallskip We set $\,U'  = X'\setminus \pi^{-1}(Z)$. 
 We denote by $E_i, i=1,\ldots ,s,$ the irreducible components of the exceptional divisor $\pi^{-1}(Q)$, and set $\, E_i^0 = E_i\cap U'.$ 
 
  We denote by $(\sE',\nabla')$ the inverse image of $(\sE,\nabla)$ on $ U' $. 
   
    \begin{prop}\label{p1} $(\sE',\nabla')$ is regular along every $E_i$.
   \end{prop}

 \smallskip  Via \ref{c2}, it is clear that \ref{t2} follows from this assertion, whose proof will occupy the next three sections.

\section{Poincar\'e-Katz rank and the Turrittin-Levelt decomposition.}\label{s3}
\subsection{} 
 Let $K$ be a field of characteristic $0$, and let $ M_K$ be a differential module over $K(({x}))$, \ie a vector space over $K(({x}))$ of finite dimension $\mu$ endowed with an action $\nabla({\theta_x})$ of ${\theta_x} = {x}\frac{d}{d{x}}$ satisfying the Leibniz rule.

\begin{ex}\label{ex1} $K = \C(E_i^0)\cong \C({y})$, ${x}$ is a local coordinate on $X'$ such that ${x}=0$ defines $E_i^0$, $y$ is a local coordinate on $E_i$, and $  M_{K }$ is the generic fiber of $ \sE' $ tensored with $\C(E_i^0)(({x}))$ over the subfield $\C(X')$, endowed with the induced action of $\nabla({\theta_x})$. \end{ex} 

According to the Turrittin-Levelt theorem, there is a finite extension $K'/K$ and a positive integer $e$ such that, putting ${x'}= {x}^{1/e}$, one has 
 a canonical decomposition of differential modules
 \begin{equation}\label{e1} M_{K }\otimes {K'(({x'}))} = \oplus  \, L_{\phi_j} \otimes_{K'(({x'}))} \, R_j \end{equation} 
 where 
\\ - $R_j$ is regular, \ie admits a basis in which the matrix of $\nabla({\theta_x})$ has entries in $K'[[{x'}]]$,
\\ - $L_{\phi_j}= K'(({x'}))$ with $\nabla({\theta_x})(1)= \phi_j\in {\frac{1}{{x'}}}K'[{\frac{1}{{x'}}}]$ (and the $\phi_j$'s are distinct).

 \smallskip\noindent  We denote by $\mu_j$ the dimension of the differential module $R_j$. 

 \smallskip Moreover, if one gathers together the summands $ \, L_{\phi_j} \otimes R_j $ according to the {\it slope}, \ie the negative of the (fractional) degree of $ \phi_j$, the resulting coarser decomposition descends to $K(({x}))$ and gives the {\it slope decomposition} 
 \begin{equation}\label{e2} M_K = \oplus_{\sigma\in \Q}\, M_{K,(\sigma)}  .\end{equation} 
 We denote by $\mu_{(\sigma)}$ the dimension of the differential module $M_{K,(\sigma)}$, so that $\mu = \sum \mu_{(\sigma)}$.

\subsection{}\label{int}  The highest slope is the {\it Poincar\'e-Katz rank of $M$}, denoted by $\rho(M)$ or $\rho$ for short.  
 If $\phi_{j}$ is of slope $\rho $, we write 
it in the form
$$ \phi_j = \phi_{j,-\rho}. {x}^{-\rho}+ h.o.t. ,\;\;  \phi_{j,-\rho}\in K'.$$
If $\rho > 0$, the coefficients $\phi_{j,-\rho}$ may be calculated as follows:
one takes a cyclic vector  $m$ for $M_K$. Then the matrix of $\nabla({\theta_x})$ in the  basis 
$(m, x^{\rho}\nabla({\theta_x})(m),\ldots, x^{(\mu-1)\rho} \nabla({\theta_x})^{\mu-1}(m))$ of $M_{K }\otimes {K'(({x'}))} $ is of the form $x^{-\rho}H$ where $H(x)\in M_\mu(K[[{x'}]])$ and the $\phi_{j,-\rho}$'s are the non-zero eigenvalues of $H(0)$ (\cf e.g. \cite[\S2]{A}). 

\subsection{  } Let us now assume that $K$ is the function field of a smooth affine curve $\Spec A$, and that $M_K$ comes from a differential module $M$ over $A(({x}))$ (as in example \ref{ex1}, where $A= \sO(E_i^{0 })$ is a localization of $ \C[{y}]$). Let $A'$ be the normalization of $A$ in $K'$.

It is not true in general that in the above decompositions, one may replace $K$ and $K'$ by $A$ and $A'$ respectively (this is the well-known problem of {\it turning points}). However, this becomes true if one restricts to a suitable dense open subset of the curve (see \cite[\S 3]{A} for a detailed analysis of this point). 

More precisely, there exists $f\in A\setminus \{0\}$ such that $A'[\frac{1}{f}]$ is finite etale over $A[\frac{1}{f}],\,\phi_j\in A'[\frac{1}{f}](({x'})) $, and the decompositions \eqref{e1} and \eqref{e2} come, respectively, from decompositions
 \begin{equation}\label{e3} M\otimes A'[\frac{1}{f}](({x'})) = \oplus  \, L_{\phi_j} \otimes_{A'[\frac{1}{f}](({x'}))}\, R_j \end{equation} 
 \begin{equation}\label{e4} M\otimes A[\frac{1}{f}](({x}))  = \oplus_{\sigma\in \Q}\, M_{(\sigma)}  .\end{equation} 
 
 \subsection{}  In general,  $\phi_j\notin A'(({x'})) $, but one always has 
 $ \phi_{j,-\rho}\in A'$, due to the above interpretation of these coefficients (\cf \ref{int}).
 
 Denoting by $\mu_j$ the rank of $R_j$\footnote{not to be confused with the dimensions $\mu_{(\sigma)}$ introduced in 3.1. Notice that $\mu_{(\rho)}= \sum \mu_j$ where the sums runs over all $j$'s such that $L_{\varphi_j}$ is of slope $\rho$, the Poincar\'e-Katz rank of $M$.}, a simple Galois argument then shows that
 $\varphi( {x}): =  \prod  ({x}^\rho - \phi_{j,-\rho})^{\mu_j}$ lies in $A[{x}]$. For $\rho>0$, this allows to define the effective divisor   $\, D  =  (\varphi( {x}))\,$  on $\Spec A[{x}]$, which is finite of degree $\mu_{(\rho)}.  \rho  $ over $\Spec A$. If $\rho=0$, we set $D=0$.
 
 \subsection{  }  When $M$ comes from an {\it integrable} connection, to the effect that there is an action $\nabla({\theta_y})$ of ${\theta_y} = {y}\frac{d}{d{y}}$ commuting with $\nabla({\theta_x})$, all the above decompositions are automatically stable under $\nabla({\theta_y})$, \ie are decompositions of integrable connections.
  
\begin{ex}\label{ex2}  
  Let us come back to example \ref{ex1}, and take $M= \sE \otimes \sO(E_i^{0 })(({x}))$. 
   In the sequel, we shall have to deal with all the $E_i$'s simultaneously, so we emphasize the index $i$ and write $\rho_i$ for the Poincar\'e-Katz rank of $\sE $ along $E_i$, and $\mu_{(\rho_i)},\, \phi_{i,j },\, \mu_{i,j},\, \phi_{i,j, -\rho_i},\, D_i, $ instead of $\mu_{(\rho )},\, \phi_{ j },\, \mu_{ j}, \, \phi_{ j, -\rho },\, D $. 
  
Geometrically, for $\rho_i>0$, the divisor 
 \begin{equation}\label{e5} D_i = ( \prod  ({x}^{\rho_i}- \phi_{i,j,-\rho_i})^{\mu_{i,j}}) \end{equation}
 has an intrinsic meaning as a {\it divisor on the normal bundle $N_{E_i^0}X'\cong \Spec A[{x}]$} (the point is that the term of lower degree of $\phi_{i,j }$ has an intrinsic meaning independent of the choice of the transversal derivation at $E_j^0$ by means of which one identifies $N_{E_i^0}X'$ and $\Spec A[{x}]$). It is {\it finite of degree $\mu_{(\rho_i)}.  \rho_i  $ over $E_i^0$}.
 \end{ex}

\section{Formal decompositions at crossings.}\label{s4}
\subsection{}\label{fd} Let us now consider a crossing point $Q'$ of $E_i$ with another component  of $\pi^{-1}(Z)$ (which is either another component $E_j$ or one of the $Z_j$'s). Let $x,y$ be etale coordinates at $Q'$ such that $E_i$ is defined by $x=0$ and the other component by $y=0$. 

The formalization of $(\sE',\nabla')$ at $Q'$ then gives rise to a differential module $\tilde M$ over $\C[[x,y]][\frac{1}{xy}]$ (with mutually commuting actions of $\theta_x$ and $\theta_y$).

\smallskip We shall say that $(\sE',\nabla')$ has {\it nice formal structure at $Q'$}\footnote{this is equivalent to saying that $Q'$ is semi-stable for $(\sE',\nabla')$ in the sense of \cite{A}, but weaker than saying that $(\sE',\nabla')$ has a good formal structure in the sense of \cite{S}.} if, after ramification along $xy=0$, there is a decomposition
\begin{equation}\label{e6} \tilde M= \oplus \, \tilde L_{\tilde\phi_k,\tilde\psi_k}\otimes \tilde R_k \end{equation}
 where 
\\ - $\tilde R_j$ is regular, \ie admits a basis in which the matrix of $\tilde\nabla({\theta_x})$ and $\tilde\nabla({\theta_y})$ have entries in $\C[[x,y]]$,
\\ - $\tilde L_{\tilde\phi_k,\tilde\psi_k}= \C[[x,y]][\frac{1}{xy}]$ with $$\tilde\nabla({\theta_x})(1)= \tilde\phi_k, \,\tilde\nabla({\theta_y})(1)= \tilde\psi_k \,\in \,\C[[x,y]][\frac{1}{xy}],\;\; \theta_x(\tilde\psi_k)=\theta_y(\tilde\phi_k).\ $$ 
Without loss of generality, one may then assume that the pairs $(\tilde\phi_k,\tilde\psi_k)$ are distinct modulo $\C[[x,y]]$ (which ensures unicity of the decomposition).

\subsection{  } For $K= \C((y))$, let us compare a decomposition \eqref{e6} with the Turritin-Levelt decomposition of the differential module $$M_K := \tilde M\otimes_{\C[[x,y]][\frac{1}{xy}]}\, K((x)).$$
One has a canonical decomposition 
 $$\frac{\C[[x,y]][\frac{1}{xy}]}{\C[[x,y]]}\;\cong \; \frac{1}{x}\C[[y ]][{1\over x}]\;\oplus \;{{{{1}\over {xy}}  \C[ {1\over x},{1\over y}]} }\;\oplus\; {1\over y }\C[[x]][{1\over y }] $$
and the projection onto the first two terms can be written 
$$\varpi :\; \frac{\C[[x,y]][\frac{1}{xy}]}{\C[[x,y]]}\to  \frac{1}{x}K[{1\over x}].$$  
  The decomposition \eqref{e6} of $\tilde M$ thus gives rise, by tensorisation with $K((x))$, to a decomposition of $M_K$ which refines \eqref{e1} (taking into account the unicity of the latter): 
 \begin{equation}\label{e6'}  L_{ \phi_j} \otimes R_j = \bigoplus_{k,\,\varpi (\tilde \phi_k)=\phi_j}  \, (\tilde L_{\tilde\phi_k,\tilde\psi_k}\otimes \tilde R_k)\otimes K((x)) .\end{equation}

 \subsection{  } Whereas it is not always the case that $(\sE',\nabla')$ has nice formal structure at $Q'$,  this can be fixed by blowing up:

\begin{prop}\label{p3} After blowing up finitely many times some crossing points, $(\sE',\nabla')$ acquires a nice formal structure at every crossing point of the inverse image of $Z$. 
\end{prop}

This was first proved by C. Sabbah \cite[III, 4.3.1]{S}, using a generalization in dimension 2 of the nilpotent orbit method of Babbitt-Varadarajan. In \cite[5.4.1]{A}, we have given a simpler and more straightforward proof (which is nevertheless too long to be repeated here).

\section{Proof of \ref{p1}.}\label{s5}  We assume from now on, as we may by \ref{p3}, that {\it $(\sE',\nabla')$ has a nice formal structure at all crossing points of  $\pi^{-1}(Z)$ lying on the exceptional divisor}. 

\subsection{  } Let $P(N_{E_i}X') $ be the projectivization of the normal bundle $N_{E_i}X'$, and let $(\infty)$ be the section at infinity over $E_i$. Taking  Zariski closure, the effective divisor $D_i$ on $N_{E_i^0}X'$ gives rise to an effective divisor $\bar D_i$ on $P(N_{E_i}X') $, which is finite over $E_i$ (of degree $\mu_{(\rho_i)}. \rho_i$). It is clear from equation \eqref{e5} that $\bar D_i$ does not meet  $(\infty)$  above $E_i^0$. 

\begin{lemma}\label{l5} $i)$ Above $E_i \cap Z'_j$, $\bar D_i$ does not meet  $(\infty)$.
\\ $ii)$ Above $Q'= E_i\cap E_{i'}$,  the intersection multiplicity of $\, \bar D_i\,$ and $\, (\infty)  \,$ is  $ \leq \mu_{(\rho_i)}.\rho_{i'} $.   
\end{lemma} 

\prf Let us use coordinates $x,y$ as in \ref{fd}. Replacing them by $x^{1/e}, y^{1/e}$ if necessary, one may assume that one has a decomposition \eqref{e6}. One reads on equation \eqref{e5} that the intersection mutiplicity of $\, \bar D_i\,$ and $\, (\infty)  \,$ above $Q'$ is 
$$  \max (0, -\sum_{j,\, \ord_x\, \phi_j = -\rho_i}\, \mu_{i,j}. \ord_y\, \phi_j).$$
It thus suffices to show that $-\ord_y\, \phi_j \leq \rho_{i'}$.
Coming back to \ref{e6}, it is clear that $$-\ord_x \tilde  \phi_k \leq  \rho_i, \; -\ord_y \tilde  \psi_k \leq  \rho_{i'},$$ and by the integrability condition  $\theta_x(\tilde\psi_k)=\theta_y(\tilde\phi_k)$, one gets 
$$-\ord_y \tilde  \phi_k \leq  \rho_{i'}, \; -\ord_x \tilde  \psi_k \leq  \rho_{i}, $$ whence (via $\varpi $) $-\ord_y\, \phi_j \leq \rho_{i'}$ (which is $0$ in case $i)$ of the lemma, since $(\sE, \nabla)$ is assumed to be regular along $Z_{i'}$). \qed

 \begin{lemma}\label{l5} For every $i$,  \begin{equation}\label{e7} \, (E_i,E_i).\rho_i\, \geq \,-\sum_{j, \,E_j\cap E_i\neq \emptyset}\, \rho_j  \end{equation}
\end{lemma}

\prf This follows from the general formula
    $$   \delta \, \deg(N_{E_i} {{X}'}) \, = \, \bar C.(0) - \bar C.(\infty)   $$ which holds for any curve (or $1$-cycle) $\bar C\subset P(N_{E_i } {{X}'})$ whose projection to $E_i$ is finite of degree $\delta $. One applies this to $\bar C = \bar D_i$, taking into account the previous lemma, noting that $\delta = \mu_{(\rho_i)}. \rho_i$ in that case, and that $ deg(N_{E_i} {{X}'})$ is the intersection number $(E_i, E_i) \leq 0$.
\qed 

\subsection{  } To finish, let us indicate how \ref{l5} implies  \ref{p1}. Let $A$ be the matrix with entries $A_{ij}= (E_i,E_j)$. We may rewrite inequality \eqref{e7} in the form
\begin{equation}\label{e8} \, \sum_j \;  A_{ij}.\rho_j \, \geq 0. \end{equation}
Now, it is well-known  that $A$ is a negative definitive symmetric matrix. Since $\rho_j \geq 0$, this together with \eqref{e8} imply that $\rho_j=0$, \ie $(\sE',\nabla')$ is regular along $E_j$.
 \qed

\section{Exponents.}\label{s6}
\subsection{} Let us come back to the situation of \ref{t2}, assuming $X$ smooth. There is a maximal open subset $V\subset X$ containing $U$ such that the complement of $U$ in $V$ is a smooth divisor (whose components are open subsets of the $Z_j$'s).   

There is a purely algebraic construction of a logarithmic lattice, \ie a locally free extension $\sE_{\log}$ of  
$\sE$ from $U$ to $V$ endowed with a logarithmic connection $\nabla_{\log}$ extending $\nabla$ with poles at $Z\cap V$,  \cf   \cite[I]{AB}. 

\begin{rem} Let $j: V\inj X$ be the open embedding. Since the complement of $V$ in $X$ is of codimension $\geq 2$, $j_\ast\sE_{\log}$ is a coherent reflexive module on $X$, hence locally free if $X$ is a surface. The difficulty of 1.3.1 is related to the fact that around $P$,  $h^\ast \nabla_{\log} $ is not a logarithmic connection on $h^\ast (j_\ast \sE_{\log})$ in general (there is a counterexample by J. Bernstein in the case where $Z$ is the union of 3 lines in the plane meeting at the origin $Q$, \cf \cite{B}). \end{rem}

 The {\it exponents} of $\nabla$ at $Z_j$ are the eigenvalues of the residues of $\nabla_{\log}$ at $Z_j\cap V$ modulo $\Z$ (the images of these eigenvalues in $\C/\Z$ do not depend on the logarithmic lattice). The following theorem is due to M. Kashiwara.

\begin{thm}\label{th2}  \cite{K}\cite{L}  The exponents modulo $\Q$ of $h^\ast \nabla$ at $P$ belong to the $\Q$-subspace of $\C/\Q$ generated by the exponents of $\nabla $ at the $Z_j$'s. \footnote{again, this statement appears in \cite[I.6.5]{AB}, but the argument given there works only in case the polar divisor has normal crossings, for the same reason as in footnote 1. } 
\end{thm}

This is easily seen to be true if the polar divisor has normal crossings, but this is more difficult beyond this case. One can reduce via \ref{l1} to the case where $X$ is a surface, or even an affine neighborhood of the origin in the complex plane. Actually, \cite{K} deals only with the case of rational exponents, but a straightforward modification of Gabber's proof in \cite{L} gives \ref{th2} as stated. 

Here is a brief sketch of Gabber's argument, coming back to the notation of 2.2. One fixes a $\Q$-linear map $\kappa: \C/ \Q\to \Q$ which sends the exponents of $\nabla $ at the $Z_j$'s to $0$. One has to show that $\kappa$ sends the exponents of $\nabla'$ at the $E_k$'s to $0$. 
 
  Let $U_k$ be a tubular neighborhood of $E_k$ in $X'(\C)$ and set $U^0_k= U_i\cap U'(\C)$. 
Then $\pi_1(U^0_k)$ is generated by elements $\alpha_i$ (turning around $E_i,\,i=1,\ldots , s$) and elements $\beta_j$ (turning around $Z'_j,\, j=1,\ldots, t$), such that $\alpha_k$ is central and $$\gamma_k := \alpha_1^{(E_1.E_k)}\ldots \alpha_i^{(E_s.E_k)} \beta_1^{(Z'_1.E_k)}\ldots \beta_i^{(Z'_t.E_k)} $$ is a commutator.  
 The monodromy representation of $\pi_1(U^0_k)$ attached to $({\sE}',\nabla')$ admits subrepresentations $V_k$ and $v_k$ where the function
 $$\chi := \kappa \circ ({\frac{1}{2\pi i}}\log)$$ applied to the eigenvalues of $\alpha_k$ takes its maximal value $M_k$ (\resp minimal value $m_k$). One has to show that $M_k= m_k=0$. 
 One has $$ \chi(det(\gamma_{k \mid V_k}))= \chi(det(\gamma_{k \mid v_k}))=0$$ (since $\gamma_k$ is a commutator) and 
 $$ \chi(det(\beta_{j \mid V_k}))= \chi(det(\beta_{j \mid v_k}))=0 $$ (by definition of $\chi$). Using the definition of $(M_k, m_k)$, this gives 
  $$\sum_j A_{ij} M_j\geq 0,\;\; \sum_j A_{ij} m_j\leq 0.$$ Since the intersection matrix $A$  is negative definite, one concludes that $M_j=m_j= 0$ for every $j$.

\begin{rem} It would be interesting to give an algebraic version of this argument in the spirit of the previous section, using residues instead of monodromy.    \end{rem}

\bigskip\bigskip \noindent {\smit Acknowledgements} {\sm. We are grateful to F. Baldassarri for many conversations about the problem treated in this paper, to J. Bernstein who prompted them by his above-mentioned counterexample, and to the referee for pointing out theorem \ref{th2} and various relevant comments. Thanks also to the organizers of the Kyoto conference for their kind invitation.}

\end{sloppypar} 

 \end{document}